%
%
\newif\ifmakeref 
\makereffalse    
%
%
%
\baselineskip = 15.5pt plus .5pt minus .5pt   
%
\font\bigbf=cmbx10   scaled \magstep1     
\font\smallcaps = cmcsc10                 
%
\ifx\fonts\cmfonts
\font\ninerm=cmr9
\font\ninei=cmmi9
\font\ninesy=cmsy9
\font\ninebf=cmbx9
\font\ninett=cmtt9
\font\nineit=cmti9
\font\ninesl=cmsl9
\else
\font\ninerm=amr9
\font\ninei=ammi9
\font\ninesy=amsy9
\font\ninebf=ambx9
\font\ninett=amtt9
\font\nineit=amti9
\font\ninesl=amsl9
\fi
\skewchar\ninei='177
\skewchar\ninesy='60
\skewchar\ninett=-1
\newskip\tglue
\def\ninepoint{\def\rm{\fam0\ninerm}
       \textfont0=\ninerm
       \textfont1=\ninei
       \textfont2=\ninesy
       \textfont\itfam=\nineit \def\it{\fam\itfam\nineit}%
       \textfont\slfam=\ninesl \def\sl{\fam\slfam\ninesl}%
       \textfont\ttfam=\ninett \def\tt{\fam\ttfam\ninett}%
       \textfont\bffam=\ninebf \def\bf{\fam\bffam\ninebf}%
       \tt\tglue=.5em plus.25em minus .15em
       \normalbaselineskip=11pt
       \setbox\strutbox=\hbox{\vrule height8pt depth3pt width0pt}%
       \let\sc=\sevenrm \let\big=\ninebig \normalbaselines\rm}%
\input amssym
\input xypic
\def\REFERENCES{locales.ref}
\input \REFERENCES
\def\fonts{cmfonts}
\advance\tolerance by 100

%
\def\runningtitlestring{The Locally Fine Coreflection and Locales}
%
\newcount\footnotecount
\footnotecount = 0
\def\makeftn[[#1]]{~\hskip-.3em{\ninepoint%
\global\advance\footnotecount by 1
\footnote{\hskip-.3em$^{\number\footnotecount}$}{#1}}}
\newcount\sectioncount
\sectioncount = 0
\def\ThisSection{\number\sectioncount}
\def\section#1{\vskip0pt plus .1\vsize
    \penalty-250\vskip0pt plus-.1\vsize\bigskip
    \global\advance\sectioncount by 1
    \noindent{\bf \number\sectioncount. #1}\nobreak\message{#1}}
\def\SectionBreak{\vskip0pt plus .1\vsize
    \penalty-250\vskip0pt plus-.1\vsize\null\bigskip}

\def\abstract#1{{\ninepoint\bigskip\centerline{\hbox{
       \vbox{\hsize=4.85truein{\noindent ABSTRACT.\enspace{#1}}
            }}}}}
\def\keywords#1{{\ninepoint\bigskip\centerline{\hbox{
       \vbox{\hsize=5.75truein{\noindent KEYWORDS.\enspace#1}
            }}}}}
\def\AmsClassMM#1{{\ninepoint\bigskip\centerline{\hbox{
       \vbox{\hsize=5.75truein{\noindent
       AMS 2000 Subject Classification Numbers:\enspace#1}
            }}}}}
\def\th #1 #2: #3\par{\medbreak{\bf#1 #2:
\enspace}{\sl#3\par}\par\medbreak}
\def\co #1 #2: #3\par{\medbreak{\bf#1 #2:
\enspace}{\sl#3\par}\par\medbreak}
\def\le #1 #2: #3\par{\medbreak{\bf #1 #2:
\enspace}{\sl #3\par}\par\medbreak}
\def\rem #1 #2. #3\par{\medbreak{\bf #1 #2.
\enspace}{#3}\par\medbreak}
\def\proof{{\bf Proof}.\enspace}

\def\sqr#1#2{{\vcenter{\hrule height.#2pt
      \hbox{\vrule width.#2pt height#1pt \kern#1pt
         \vrule width.#2pt}
       \hrule height.#2pt}}}

\overfullrule=10pt
\def\boxit#1{\vbox{\hrule \hbox{\vrule \kern2pt
                 \vbox{\kern2pt#1\kern2pt}\kern2pt\vrule}\hrule}}
%
\newdimen\refindent\newdimen\plusindent
\newdimen\refskip\newdimen\tempindent
\newdimen\extraindent
\newcount\refcount
\newwrite\reffile
\def\beginref{\ifmakeref\immediate\openout\reffile = \REFERENCES\else\fi}
\def\endref{\ifmakeref\immediate\closeout\reffile\else\fi}
\def\ifundefined#1{\expandafter\ifx\csname#1\endcsname\relax}
\def\referto[#1]{\ifundefined{#1}[?]\else[\csname#1\endcsname]\fi}
%
%
\refcount=0
\def\ref#1:#2.-#3[#4]{\ninepoint 
\advance\refcount by 1
\setbox0=\hbox{[\number\refcount]}\refindent=\wd0
\plusindent=\refskip\extraindent=\refskip
\advance\plusindent by -\refindent\tempindent=\parindent %
\parindent=0pt\par\hangindent\extraindent %
 [\number\refcount]\hskip\plusindent #1:{\sl#2},#3
\parindent=\tempindent
\ifmakeref
\immediate\write\reffile{\string\expandafter\def\noexpand\csname#4\endcsname%
{\number\refcount}}\else\fi}
\refskip=\parindent
%

\def\makeknown#1#2{\expandafter\gdef
            \csname#2\endcsname{\hbox{\csname#1\endcsname #2}}}
\makeknown{rm}{Cov}
\makeknown{rm}{int}
\makeknown{rm}{Ord}
\makeknown{rm}{Rank}
\makeknown{rm}{rank}
\makeknown{rm}{Root}
\makeknown{rm}{End}
\makeknown{rm}{Level}
\makeknown{rm}{Seg}
\makeknown{rm}{pr}
\makeknown{rm}{top}
\makeknown{rm}{cf}
\makeknown{rm}{loc}
\makeknown{rm}{Hom}
\makeknown{rm}{Pt}
\makeknown{bf}{Loc}
\makeknown{rm}{op}
\makeknown{bf}{Top}
\makeknown{bf}{p}
\makeknown{sans}{COV}
\def\c#1{{\cal #1}}

\def\Cech{\check{\cal C}}

\def\twolattice{\hbox{\bf 2\hskip-.4em}|} 
%

\def\specialheadlines{
        \headline={\vbox{\line{
                     \ifnum\pageno<2 \ffolio
                     \else\rightheadline
                     \fi}\bigskip
                     }}}
        \def\ffolio{\hfil}
        \def\rightheadline{\hfil
                    {\smallcaps\runningtitlestring}\hfil
                    \hbox{\rm\folio}}
\nopagenumbers
\specialheadlines
%
%
\null
\vskip1.5truecm
  \centerline{\bigbf On the Locally Fine Construction in Uniform Spaces,}
  \medskip
  \centerline{\bigbf Locales and Formal Spaces}
\vskip.75truecm
 \centerline{\sl Aarno Hohti, University of Helsinki}
\vskip.5truecm
\vskip.5truecm
\abstract{We investigate the connection between 
          the spatiality of locale products and
          the earlier studies of the author on the locally fine
          coreflection of the products of uniform spaces. After
          giving a historical introduction and indicating the
          connection between spatiality and the locally fine
          construction, we indicate how the earlier results directly solve
          the first of the two open problems 
          announced in the thesis of T.\ Plewe. Finally,
          we establish
          a general isomorphism between the
          covering monoids of the localic product of topological
          (completely regular) spaces and the locally
          fine coreflection of the corresponding product of
          (fine) uniform spaces. Additionally, paper relates
          the recent studies on formal topology and uniform
          spaces by showing how the transitivity of covering
          relations corresponds to the locally fine construction.}

\smallskip
\AmsClassMM{Primary 54E15, Secondary 54B10.}
\vskip.5truecm
%
\section{Introduction.} This paper\hskip-.25em\makeftn[[During 
the initial research for this paper, the author
was assisted by a grant from The Finnish Academy of Sciences and
the grant no.\ 201/97/0216 from the Grant Agency of the Czech Republic.]]
is based on the work carried out on the products of uniform spaces 
since 1981 by the author and others (see \referto[Hohtia] -- \referto[Hohtig],
\referto[HohtiPelanta], \referto[HohtiPelantb], \referto[HohtiYun]) using
the technique of the locally fine corecflection. We 
study the connection of these results with the spatiality of
localic products of topological spaces and the relationship of 
the locally fine corecflection with the so-called formal spaces.

Our work has been
motivated by the appearance
in 1996 of T.\ Plewe's Ph.\ D.\ thesis 
{\referto[Plewea]\makeftn[[The 
doctoral dissertation, on which \referto[Plewea]
is based, had appeared in 1994.]]} on the
spatiality of localic products of topological, especially separable
metrizable spaces. 
Not only did the author discover an analogy between
the previous results and those obtained in the thesis in
question, but the solution of the first open problem left open
(re-solved consistently with ZFC) in \referto[Plewea] was
seen to follow from the author's article \referto[Hohtif].
Therefore, we decided to write this paper in order to
indicate the existence and usefulness of such previous
research in an equivalent field. Let us note that the
link between the spatiality of the localic products of
paracompact spaces and the corresponding product
of fine (uniform) spaces was already established by J.\ Isbell in
\referto[Isbellc].

On the other hand, the connection between locales and
so-called {\it formal spaces} as clearly indicated in
the thesis of I.\ Sigstam(\referto[Sigstama]), led the
author to discover the relation of the latter to
(pre-)uniform spaces defined by means of filters of
coverings of a space. The transitivity of the
covering relations of formal spaces is obtained by
applying the locally fine construction. Hence, we obtain
an interesting link between three `unorthodox' approaches
to topology, viz.\ locales (frames), uniformities and
formal spaces. In fact, the locally fine coreflection
was first introduced and studied by Ginsburg and Isbell
in the beginning of the 1950's as a combinatorial
approach to topology (given a monoid or filter of
coverings). This project failed as one obtained topology
-- for metric spaces -- only in the case of complete spaces,
or more generally, when the spaces considered were
paracompact, in the case of so-called {\it supercomplete}
spaces. Their work \referto[Ginsburga] was published
in 1959.

However, the then recently introduced study of `local lattices'
and `paratopologies' (see J.\ Benabou \referto[Benaboua],
D.\ and S.\ Papert \referto[Paperta]) was taken up by
Isbell resulting in an article on uniform locales
\referto[Isbellc], in which he characterized \referto[Ginsburga]
as a paper ``about objects in a category $\c{H}$ now visible
as the hypercomplete uniform locales'' (ibid., p.\ 31). The hypercomplete
uniform locales obtained from uniform spaces also have
spatial products, provided that the corresponding topological products
are paracompact. The study of locales (their opposite objects
being named frames) is now a well-established field,
both in topology (Johnstone) and uniform spaces, closely
related to topos theory (see MacLane and Murdijk).
Isbell's student Plewe extended results on
spatial products, and also improved our theorem in
\referto[HohtiYun] that \v Cech-scattered paracompacta have
paracompact countable products by establishing that
countable localic products of partition-complete
(\referto[Michaelb], \referto[TelgarskyWicke]) paracompact spaces are spatial.
The connection to the locally fine coreflection was,
however, never pointed out.

Formal spaces form a counterpart of the original
constructive approaches to topology which
considered recursive sequences
to define points: One gives a collection of `pieces' of a
space related through a covering relation and
studies recursive constructions on the basis of
these pieces. After being introduced by Fourman and Grayson in
1982 (\referto[FourmanGrayson]) and made manifest by
Sambin in 1987 (\referto[Sambina]), this approach
to a point-free topology has been studied by several
authors (see, e.\ g., Negri and Valentini \referto[NegriValentini]
and Sigstam \referto[Sigstama]). As will be seen in the sequel,
the transitivity requirement of the covering relation
of a formal space is the counterpart of Kuratowski's
classical condition that $C^2=C$ for a closure operation $C$,
and is essentially equivalent to the locally fine construction
$\lambda$, when the definition is extended from
pre-uniformities (covering monoids) to covering relations.
This connection will perhaps give a justification for
the original attempt by Ginsburg and Isbell to obtain
topology combinatorially through $\lambda$.

\section{The locally fine coreflection.} Open covers of 
topological spaces have an obvious local character in the 
following sense: If $\c{G}$ is an open cover of a topological space,
and for each $G\in \c{G}$, $\c{H}_G$ is an open cover of $G$,
then again the combined family $\cup\{\c{H}_G: G\in \c{G}\}$ is
an open cover of the space. This is, however, not valid when `open' is
replaced by `uniform'. The locally fine operation $\lambda$ can be
thought of as an attempt to reach topology from a given filter
of coverings through combinatorial localization, i.\ e., by
closing the given filter under the above condition.
Let $\mu$
be such a filter of coverings (pre-uniformity) on
a set $X$. Assume that $\{U_i\}$ is a member of $\mu$ and for
each $i$, we are given a member $\{V^i_j\}$. Then
the above condition requires that the `uniformly locally uniform'
cover $\{U_i\cap V^i_j\}$ be again uniform. A pair $(X,\mu)$
with this property is called {\it locally fine}. In case
$(X,\mu)$ is not locally fine, we may define
the closure of $\mu$ under this construction as follows.

For generality, let $\nu$ be another pre-uniformity on the set $X$.
Then the filter $\mu/\nu$ is defined to consist of all
coverings having a 
refinement of the form $\{U_i\cap V^i_j\}$, where ${U_i}\in\mu$
and for each $i$, ${V^i_j}\in \nu$. Now let us define
by transfinite iteration the consecutive {\it Ginsburg-Isbell}
derivatives\makeftn[[Ginsburg and 
Isbell (\referto[Ginsburga]) define $\mu^{(\alpha + 1)}$ 
as $\mu^{(\alpha)}/\mu^{(\alpha)}$ which, however, is less 
suitable for certain inductive purposes. The above 'slowed down' 
version
was introduced by the author in \referto[Hohtic]. For 
products of the form $\mu\times \nu$, we may still slow 
down the derivation by proceeding one coordinate at a time:
Instead of using covers of the form $\c{U}\times \c{V}$, 
where $\c{U}\in \mu$ and $\c{V}\in \nu$, one considers
covers $\c{U}\times \{Y\}$, $\{X\}\times \c{V}$. The corresponding 
`coordinatewise refinement condition' can be extended 
to arbitrary products.]] 
by setting $\mu^{(0)} = \mu$;
$\mu^{(\alpha + 1)} = \mu^{(\alpha)}/\mu$, and
$\mu^{(\beta)} = \cup\{\mu^{(\alpha)}: \alpha < \beta\}$ when
$\beta$ is a limit ordinal. There will be a least $\alpha$ such that 
$\mu^{\alpha+1} = \mu^{\alpha}$; this filter is called the 
{\it locally fine coreflection} of $\mu$ and denoted 
$\lambda\mu$. On of the essential results in \referto[Ginsburga]
state that if the filter $\mu$ is a uniformity, then so is 
$\lambda\mu$. 

Is $\lambda\mu$ sufficient for the topology of the underlying 
space $X$, even in the case of a metric uniformity? There is a 
curious connection with the completeness of hyperspaces, 
hence the term {\it `hypercomplete'} or {\it `supercomplete'}. In case 
the metric uniform space $\rho X$ is complete, then by 
\referto[Ginsburga] $\lambda\rho$ is the fine uniformity 
of $X$. (The fine uniformity is the collection of all normal
covers of a given completely regular space, or in terms of 
entourages the filter of all neighbourhoods of the diagonal
$\Delta(X)$ in $X\times X$.) In this case every open cover 
is in $\lambda\rho$, because (by A.\ H.\ Stone) every open cover 
of a metrizable space is normal. It was previously known that 
the hyperspace of a complete metric space is again
complete. Isbell proved that for a uniform space $\mu X$,
the uniform hyperspace $H(\mu X)$ is complete if, and only if,
the locally fine coreflection $\lambda\mu$ is fine and $X$
is paracompact (\referto[Isbella]). Thus, for such supercomplete
spaces, the locally fine construction is sufficient for defining
topology (open covers) combinatorially from uniform ones.

\section{Locales.} For a topological space $X$, the topology of 
$X$, written $T(X)$, is a complete lattice which additionally satisfies
the following Heyting axiom:
$$ x\wedge \bigvee y_\alpha \; = \; \bigvee x\wedge y_\alpha. \leqno(*)$$
%
%
%
%
%
Given topological spaces $X,Y$ and a continuous 
mapping $f:X\to Y$, there is a 
natural homomorphism $f^*: T(Y)\to T(X)$ obtained 
by sending an open set $O\in T(X)$ to its preimage
under $f$. 
Complete lattices satisfying the axiom $(*)$ 
with opposite morphisms are called {\it local lattices}
(Benabou) or simply {\it locales} (Isbell). Thus, we have
a contravariant functor $T: \Top\to \Loc$.
Very general 
`local structures' satisfying $(*)$ were considered 
by Ehresmann in \referto[Ehresmanna], who also defined
the notion of `paratopology', studied in
the papers \referto[Paperta], \referto[Benaboua].
The category of {\it frames}
is the opposite $\Loc^{\op}$ of the category 
of locales, i.\ e., they are complete right distributive 
lattices with morphisms in the standard direction.


The other way goes from locales to spaces.  
Recall that the {\it points} of a locale $L$ may be considered
homomorphisms $\phi: L\to \twolattice$, where $\twolattice$
denotes the lattice with two elements $0,1$ such that 
$0 < 1$.
One defines the topological 
space $\Pt(L)$ of points by choosing for the subbasis 
the sets 
$ x^* \;=\; \{\phi\in \Hom(L,\twolattice): 
                     \phi(x) = 1\}.$
A question arises concerning the relations of 
$X$, $T(X)$ and $\Pt(T(X))$. One says that $L$ is 
{\it spatial} (has enough points) if distinct 
elements of $L$ can be distinguished by points,
i.\ e., if for $x,y\in L$, $x\neq y$, there is 
a point $\phi: L\to \twolattice$ such that
$\phi(x)\neq\phi(y)$. Then if $L$ has enough
points, $L$ is isomorphic to $T(X)$ for some 
space $X$. Isbell has proved e.\ g.\ that 
quasi-compact regular (more generally `subfit')
locales are spatial
(\referto[Isbellc], 2.1). On the other hand, 
spatiality is hard to preserve in the products 
of locales. We postpone the constructive 
definition of localic products to the last section
where we prove our main result. We simply note here
that the category of locales has products, 
denoted here with the symbol $\Pi_{\loc}$. (Similarly,
the category of frames has co-products.) One always
has the equality (topological homeomorphism)
$$ \Pt(\Pi_{\loc} L_i) \; \cong \; \Pi \Pt(L_i).$$
Even assuming that the $L_i$ come from spaces, i.\ e., 
$L_i = T(X_i)$ and the spaces $X_i$ are `sober' so that
$\Pt(L_i)$ is homeomorphic to $X_i$, this is still
not enough. We can only infer that 
$$ \Pt(\Pi_{\loc}T(X_i)) \; \cong \; \Pi X_i,$$
and this does not say anything about the {\it pointless part}
of the localic {\it product}. We need to establish that 
the spaces $X_i$ are preserved under the localic product,
i.\ e., $\Pi_{\loc} T(X_i) \cong T(\Pi X_i)$. We often 
replace such an isomorphism with an equality.

The following section seeks to show the relation of 
the studies (mainly in the 1980's) on the locally
fine coreflection of product uniformities to the
thesis of T.\ Plewe.

\section{Product theorems.}
Detailed studies on the behaviour of $\lambda$ on complete metric
and other uniform spaces were carried out by the author and 
J.\ Pelant in the 1980's. Let us first consider supercomplete 
spaces. The questions considered were mainly related to product
spaces, a topic directly connected with (and in the metrizable case
equivalent to) the spatiality of localic products (Plewe's thesis
\referto[Plewea]). If the topological spaces $X,Y$ are paracompact
(and Hausdorff), then the fine uniform spaces $\c{F}X, \c{F}(Y)$
are supercomplete. The condition for the product to be supercomplete
is that it is (topologically) paracompact and the equation
$$ \lambda(\c{F}(X)\times \c{F}(Y)) \; = \; \c{F}(X\times Y)) $$
holds. In the infinite case, the corresponding equation has the 
form $\lambda\Pi\c{F}(X_i) = \c{F}(\Pi X_i)$. The main results
proved in the series {\it On Supercomplete Spaces I -- V}  
are listed below. A space $X$ is called $\c{K}$-scattered 
with respect to a class $\c{K}$ of spaces if every non-empty
closed subspace of $X$ contains a point with a closed neighbourhood which 
belongs to $\c{K}$. In the sequel, 
$\c{C}$ (resp.\ $\Cech$) denotes the class of compact 
(resp.\ \v Cech-complete) spaces.

\th Theorem 1. (\referto[Hohtid], \referto[Hohtig]): A binary product 
  $\c{F}(X)\times \c{F}(Y)$ is 
  supercomplete for every paracompact space $Y$ if, and only if,
  $X$ is paracompact and $\c{C}$-scattered. 

  This characterization
  was partially obtained in \referto[Hohtid], and completed 
  in \referto[Hohtig]. The paper \referto[Hohtid] (\referto[HohtiReport])
  contained the result that $\c{F}(X)\times \c{F}(Y)$ is 
  supercomplete for every paracompact $Y$ whenever 
  $X$ is C-scattered and paracompact. Furthermore, a 
  partial converse obtained stated that if $X$ is a \
  paracompact $p$-space (of Arhangel'skii), then 
  $\c{F}(X)\times \c{F}(S)$ is supercomplete for 
  every separable metrizable $S$ iff $X$ is C-scattered. 
  Therefore, a metrizable space $X$ is a {\it multiplier}
  in the class of supercomplete (topologically) metrizable spaces 
  iff it is $C$-scattered. This result implies 
  one of the results obtained in Plewe \referto[Plewea], namely
  that for metrizable spaces, the locale $X\times_{\loc} Y$
  is spatial for all $Y$ if and only if $X$ is completely
  metrizable and does not contain a closed copy of the 
  irrationals. (See below for a discussion 
  on the spatiality of metrizable products.) 
  Indeed, for metrizable spaces, the properties of 
  1) being C-scattered and 2) being completely metrizable 
  and not containing a closed copy of the irrationals 
  are equivalent. (If $X$ is $C$-scattered, then by
  \referto[Telgarskya], Theorem 1.7, 
  it is an absolute $G_\delta$ space, and hence completely
  metrizable. The space $X$ cannot contain a closed 
  copy of of the irrationals $\Bbb J$, because $\Bbb J$
  is nowhere locally compact. On the other hand, if 
  $X$ is not $C$-scattered, then $X$ contains a closed subspace
  $F$ which is nowhere locally compact. Then by \referto[Hausdorffa], 
  p.\ 157 (or (\referto[Plewea], 2.1), if $X$ is completely
   metrizable, $F$ contains a closed copy of $\Bbb J$.)

  We note here that Isbell (\referto[Isbelld], Th.\ 4) proved 
  that in the class of {\it completely regular spaces}, 
  $X\times_{\loc} Y$ is spatial for all $Y$ if, and only if, 
  $X$ is C-scattered. 6

\th Theorem 2. (\referto[HohtiPelantb]): If $X$ is paracompact 
  and $\c{C}$-scattered, then 
  the countable power $\c{F}(X)^{\Bbb N}$ is supercomplete.
  
  Indeed, if $X$ is assumed to be
  merely a countable union of closed $\c{C}$-scattered 
  subspaces, then the so-called {\it metric-fine} coreflection 
  $m(\c{F}(X)^{\Bbb N})$ is supercomplete. 
  This gave a new topological 
  corollary for paracompact spaces, made possible by
  the Noetherian tree technique (see below) related to the operation 
  $\lambda$. Well-known topological cases known previously
  (the locally compact paracompact case in 
   Arhangels'kii \referto[Arhangelskiia] and Frol\'\i k
   \referto[Frolika], $\c{C}$-scattered 
  Lindel\"of  spaces (Alster, \referto[Alstera]), 
  scattered paracompact spaces (Rudin, Watson \referto[RudinWatson]) and 
  paracompact $\c{C}$-scattered spaces (Friedler, Martin,
  Williams \referto[Friedlera])) all follow from the uniform 
  case by suitably choosing the uniformity considered.
  The result of \referto[HohtiPelantb] was extended to
  $\Cech$-scattered paracompact spaces by the author and Yun Ziqiu
  in \referto[HohtiYun], first announced in 1990 (conferences 
  in Dubrovnik and Tsukuba). This time the topological 
  corollaries were new: The product theorem holds for 
  $\Cech$-scattered Lindel\"of, paracompact, and ultraparacompact
  spaces. 

\th Theorem 3. (\referto[HohtiYun]): If $X$ is a $\Cech$-scattered 
   paracompact space, then the countable power $\c{F}(X)^{\Bbb N}$
   is supercomplete.

On the other hand, the author proved `omitting' 
theorems for supercompleteness in products. The method was 
based on the notion of {\it $n$-cardinality}, due to 
T.\ Przymusinski and van Douwen, who gave
similar applications to topological spaces 
(cf.\ \referto[Przymusinskia]). 

\th Theorem 4. (\referto[Hohtie]): For each $n \in \{1,2,3,\ldots\}$ 
   there is a subset $X\subset \Bbb R$ such that 
   $\c{F}(X)^k$ is supercomplete for $k = 1, \ldots, n$
   but $\c{F}(X)^{n+1}$ is not, in other words 
   $\lambda(\c{F}(X)^n) = \c{F}(X^n)$ but 
   $\lambda(\c{F}(X)^{n+1}) \neq \c{F}(X^{n+1})$.

In the same paper, it was also established that the set 
$X$ can be chosen so that all finite powers of $\c{F}(X)$
are supercomplete, while the countable power is not. 
As mentioned in \referto[Hohtie],
the sets $X$ were constructed as Bernstein sets. One of the 
corollaries in Plewe's thesis (\referto[Plewea]) 
is a result that follows from the above theorem, and is 
in fact directly equivalent to it:

\th Theorem 4'. (Corollary 5.6 in \referto[Plewea]): For each 
  $n\in \{2,3,\ldots\, \omega\}$ there exist Bernstein 
  sets whose $m$th localic power is spatial for each $m<n$,
  while the $n$th localic power is not.

Indeed, the equivalence of spatiality and supercompleteness
was pointed out already by Isbell in 1972 (\referto[Isbellc],
Theorem 3.12): The product locale of supercomplete spaces $X_i$
is the locale underlying the hypercompletion (as a locale) 
of their product space. 
Thus, the product locale 
is the locale derived from the product space (and hence the 
product is spatial) if, and only if,
the product space itself is hypercomplete. As countable products
of metrizable spaces are always paracompact, the product 
of at most countably many metrizable $X_i$ is spatial
(as locale) if, and only if, 
$\lambda(\Pi \c{F}(X_i)) = \c{F}(\Pi X_i)$. 

This equivalence is not, however, the end of the story. In his 
thesis Plewe listed two unanswered questions, of which the 
first is directly related to Theorems 4--4': {\it Do there 
exist non-complete spaces with spatial countable localic powers?}  
He proved that the question has a positive answer in case 
one assumes a set-theoretical hypothesis consistent with the ZFC,
namely that $|{\Bbb R}| \geq \omega_2$ and the unions of $\omega_1$
first category subsets is again of the first category
(this is implied by Martin's Axiom). However, the author had 
extended the technique of $n$-cardinality and published a 
solution to the equivalent problem for uniform spaces 
in 1988:

\th Theorem 5. (\referto[Hohtif], 3.2): There is a non-analytic subset
  $X\subset [0,1]$ such that $\lambda(\c{F}(X)^\omega) = 
  \c{F}(X^\omega)$.

As a corollory (\referto[Hohtif], 3.3), one directly obtains 
from Gleason's factorization theorem (cf.\ \referto[Isbellc],
p.\ 130) that the equality in the above theorem is valid 
for any power. For definiteness, let us give here the 
corollary to Theorem 5 for spatiality:

\th Theorem 5': There is a non-analytic and hence non-complete
subset $X\subset [0,1]$ such that the countable localic power 
of $T(X)$ is spatial.

The original notion of $n$-cardinality was used to extend 
the validity of the CH for Borel sets (Alexandroff) to the 
$n$-cardinality version of the CH for subsets of finite 
products. Let $X$ be an arbitrary set and let $A\subset X^n$.
Consider finding a set $Y\subset X$, as small as possible,
such that the codimension 1 `hyperplanes' $\pi_i^{-1}(y), y\in Y$
cover the set $A$. Accordingly, we define the $n$-cardinality 
of $A$, written $|A|_n$ as the minimum cardinality of a subset
$Y\subset X$ such that 
$$ A\subset Y\times X^{n-1} \cup X\times Y\times X^{n-2}
    \cup \cdots \cup X^{n-1}\times Y,$$
or, equivalently, $A\subset \cup\{\pi_i^{-1}[Y]: 1\leq i\leq n\}$.
The result proved by Przymusinski in \referto[Przymusinskia]
states that if $X$ is a Polish space and $A\subset X^n$ is 
an analytic subset with $|A|_n > \omega$, then the 
$n$-cardinality equals $2^\omega$.

For dealing with infinite powers, the author defined in 
\referto[Hohtif] the notion of {\it relative} $\omega$-cardinality: 
The $\omega$-cardinality of a subset $A\subset X^\omega$
with respect to a subset $S\subset X$, written $|A,S|_\omega$,
is the minimum cardinality of a subset $Y\subset S$ 
(if such a set exists) such that 
$$ A \subset \cup \{\pi_i^{-1}[Y]: i\in \omega\}.$$
In case there is no such $Y\subset S$, we define 
$|A,S|_w = |X|$. Due to the relativity condition, this 
is a non-trivial extension of the notion of $n$-cardinality.
(On the other hand, a similar notion of relative 
$n$-cardinality permits a simple proof of the basic 
result, see \referto[Hohtif], 2.1.) The main principle in the 
inductive proof of Theorem 5 was the result that
if $X$ is a Polish space, $S\subset X$ is arbitrary, 
and $A\subset X^\omega$ is analytic, then 
$|A,S|_\omega > \omega$ implies $|A,S|_\omega = 2^\omega$.

\noindent
{\bf Remark 1:} A basic example of a non-spatial product is 
given by $\Bbb Q\times \Bbb Q$. This example was explicitly
handled by Johnstone in his book (\referto[Johnstonea], II 2.14), which 
appeared in 1982. Coincidentally, in the same year, the author had 
shown as a particular corollary to his results on supercompleteness
that $\lambda (\c{F}\Bbb Q\times \c{F}\Bbb Q) \neq
   \c{F}(\Bbb Q\times \Bbb Q).$ This followed from the 
following result: Given Tychonoff spaces $X,Y$ such that
$X\times Y$ is Lindel\"of, then 
$\lambda(\c{F}(X)\times \c{F}(Y)) = \c{F}(X\times Y)$ if,
and only if, for each compact $K\subset (\beta X\times \beta Y)
\setminus (X\times Y)$ there are \v Cech-complete 
paracompact subspaces $M,N$
of $\beta X, \beta Y$, respectively, such that 
$X\subset M, Y\subset N$ and $K\cap (M\times N) = \emptyset$.
(See \referto[Hohtid], 3.5). For a subset $X$ of 
the unit interval $I$,
there is an easier way of paraphrasing this result: 
$\c{F}X\times \c{F}X$ is supercomplete if, and only if,
for each compact $K\subset I^2\setminus X^2$
there is a first category subset $A\subset I\setminus X$ such that 
$$ K\subset (A\times I) \cup (I\times A).$$
However, there is an entire geometric circle 
$C\subset I^2\setminus {\Bbb Q}^2$, and it cannot 
be covered by the projection pre-images
$\pi^{-1}_i[A]$, $i=1,2$ of any first category set $A$.

It was also shown that $\c{F}({\Bbb J})\times \c{F}(Q)$, where 
$\Bbb J$ denotes the irrationals, is not supercomplete.
Thus, it follows that whenever a Tychonoff space $X$ 
contains a closed copy of the irrationals, then 
$\c{F}(X)\times \c{F}({\Bbb Q})$ is not supercomplete, 
and the localic product $X\times_{\loc} {\Bbb Q}$ is not
spatial.

\section{Noetherian trees.} The method used in proving 
the positive countable product theorems 2--3 was based on
trees with only finite branches. The application of the 
locally fine condition in the successive constructions 
of the covers in the derivatives $\mu^{\alpha}$ leads to 
such `Noetherian' covering trees. In such a tree, the 
immediate successors of an element form its uniform cover,
and the collection $\End(T)$ of all maximal elements of the 
tree $T$ cover the space. This technique was first used in 
Pelant's proof (\referto[Pelanta])\makeftn[[Z.\ Frol\'\i k
had an interesting interpretation of this result: The fine 
spaces (i.\ e., Tychonoff topologies) form the smallest 
coreflective class such that all subspaces are 
locally fine. Thus, Tychonoff spaces are obtained 
from the locally fine spaces through a purely 
categorical construction.]]
of Isbell's conjecture that every locally fine 
space is `subfine' (a subspace of a fine space). Pelant 
showed that the `$\lambda$ equation' considered above,
$\lambda (\Pi\c{F}(M_i)) = \c{F}(\Pi M_i)$ holds for 
{\it any} collection of completely metrizable spaces 
$M_i$. (See below for a current extension of this result.)
Noetherian trees were used to represent the recursive
construction of covers in the consequtive derivatives
$\mu^{(\alpha)}$. The essential lemma used by Pelant 
states (in our formulation) that $\c{U}\in \lambda\mu$ 
if, and only if, there
is a Noetherian tree $T$ of subsets of the underlying space 
such that 1) $T$ satisfies the {\it uniform covering 
condition} with respect to $\mu$ (i.\ e., the immediate 
successors of a non-maximal element form its
uniform cover) ; 2) the maximal 
elements $\End(T)$ form a cover which refines $\c{U}$
and 3) $T$ has $X$ as its root. Each cover 
$\c{G}\in \mu^{(\alpha)}$ can be reached by such a  
Noetherian tree and vice versa.
This enables one to replace the consecutive derivatives
and transfinite induction by arguments based on 
well-foundedness. It should be noted that general 
(localic) products of completely metrizable spaces 
(not being paracompact) are not spatial; the 
equation $\lambda (\Pi\c{F}(M_i)) = \c{F}(\Pi M_i)$ is 
not sufficient alone but must be complemented 
with the condition that each open cover of the 
product is normal. We will give a more general result
in the last section.

\noindent
{\bf Remark 2:} Noetherian trees have well-defined {\it ranks}, and 
complete metric spaces of a finite or countable 
rank were studied by the author in \referto[Hohtic]. 
(We say the rank of a complete metric space $\rho X$
is the least $\alpha$ such that $\rho^{(\alpha)} = 
\c{F}(X)$, the existence of which is quaranteed by 
\referto[Ginsburga], 4.2.) Among other results, it was proved that 
for a finite or countable $\alpha$, the rank of 
$\rho X$ equals $\alpha$ if, and only if, $X$ has a 
compact set $K$ such that outside of any neighbourhood
of $K$, $X$ is uniformly locally of a strictly lesser
rank. This naturally led the author to recursively 
constructed decompositions of such spaces into 
Noetherian trees of closed subspaces in which the 
maximal elements are compact. The extended 
results obtained in \referto[HohtiPelanta] by 
the author and Pelant have to be bypassed here.

\section{A game-theoretical characterization.}
Noetherian trees can be used to give a direct motivation
to a game-theoretical characterization of supercompleteness
introduced -- but not studied -- by the author in 1983
\referto[HohtiEger]. There are two players I and II. For 
each game we choose an open cover $\c{V}$ of the given 
uniform space $\mu X$. Player I begins by choosing a 
uniform cover $\c{U}_0\in \mu$. If possible, Player II 
responds by selecting an element $U_0\in \c{U}_0$ such that 
$U_0\subset V$ for no $V\in \c{V}$. Then Player I continues 
by choosing a uniform cover $\c{U}_1$ of $U_0$. Player II
again selects -- if possible -- an element $U_1\in \c{U}_1$
such that no member of $\c{V}$ contains this $U_1$. 
Inductively, after the choice $U_n$ by Player II, 
Player I chooses a uniform cover $\c{U}_{n+1}$ of 
$U_n$ and Player II selects, whenever possible, 
an element $U_{n+1}\in \c{U}_{n+1}$ such that 
$U_{n+1}\subset V$ for no $V\in \c{V}$. Otherwise, the 
play stops at $U_n$. If this play of the game
$G(\mu X, \c{V})$ has infinitely many moves, then Player II
wins, otherwise Player I wins. Then we may state the 
following characterization of supercompleteness in terms 
of the games $G(\mu X, \c{V})$:

\th Theorem 6. (\referto[HohtiEger], Theorem 5'): A uniform 
  space $\mu X$ is supercomplete if and only if, for any open 
  cover $\c{V}$ of $X$, Player I has a winning strategy
  in the game $G(\mu X, \c{V})$.

\proof If $\mu X$ is supercomplete, then $\c{V}\in \lambda\mu$,
 and there is a Noetherian tree $T$ with $\Root(T) = X$, 
 $T$ satisfies the uniform covering condition and 
 $\End(T) \prec \c{V}$. By proceeding along the branches of 
 $T$, and using the uniform covering condition, Player I 
 has a (stationary) winning strategy in the game $G(\mu X, \c{V})$.
 
 On the other hand, suppose that Player I always has such a 
 winning strategy. Given an open cover $\c{V}$ of $X$, 
 it is enough to produce a Noetherian tree $T$ as in the 
 preceding paragraph. As Player I has a winning strategy 
 in $G(\mu X, \c{V})$, one is able to find a uniform 
 cover $\c{U}$ of $X$ such that Player I knows how to 
 win every play following Player II selecting elements
 $U\in \c{U}$. The construction of $T$ stops at every 
 $U\in \c{U}$ which is contained in some member of $\c{V}$.
 (Those are choices that Player II cannot make.) On the other
  hand, we will continue with all other members of $\c{U}$.
 Player I chooses, for each such member a uniform cover, 
 and the definition of $T$ is inductively continued. 
 By taking the union of all the inductive steps, we get 
 a Noetherian tree $T$, because each branch corresponds 
 to a play of $G(\mu X, \c{V})$ in which Player I wins.
 By the construction of $T$, we have $\End(T) \prec \c{V}$,
 as required.
 
Thus, the winning strategy in the game $G(\mu X, \c{V})$ 
is directly obtained from the Noetherian tree associated 
with any refinemenent of $\c{V}$ in $\lambda\mu$. Each 
particular play can be won by Player I by following a 
particular branch of such a tree. A game-theoretical 
characterization of spatiality in localic products 
$X\times_{\loc} Y$ was given by Plewe in \referto[Plewea],
likewise related to trees (\referto[Plewea], p.\ 647). 

By applying the above theorem to products of uniform 
spaces, we immediately 
obtain a characterization of their supercompleteness
as follows. It is enough to consider the case in which 
the factors $X,Y$ are fine paracompact spaces. In the game 
$G(X,Y, \c{W})$, we are given an open cover $\c{W}$ of 
$X\times Y$. We may assume -- if necessary -- that $\c{W}$
consists of open rectangles $W_1\times W_2$. Player I 
chooses open covers $\c{U}_0, \c{V}_0$ of $X$ and $Y$,
respectively, claiming that the rectangular cover 
$\c{U}_0\times \c{V}_0$ refines $\c{W}$. Player II 
selects, if possible, a rectangle $U_0\times V_0$, 
$U_0\in \c{U}_0$, $V_0\in \c{V}_0$, such that 
$U\times V$ is not contained in any member of $\c{W}$.
Then Player I chooses open covers $\c{U}_1$, $\c{V}_1$
of $U_0$ and $V_0$, respectively, obtained by restricting 
open covers of $X$ and $Y$, and claims 
that $\c{U}_1\times \c{V}_1$ refines the restriction of 
$\c{W}$ to $U_0\times V_0$. The rest of the play
is defined inductively, and Player I wins, if it only 
involves finitely many moves; otherwise, Player II wins. 
Again, the product is supercomplete if Player I has 
a winning strategy in $G(X,Y,\c{V})$ for each open cover
$\c{V}$. 

It is not directly possible to change the rules of the games 
$G(X,Y,\c{V})$ so that Player I would choose simple rectangles
$U_i\times V_i$, instead of choosing covers. The crux of 
the rules is to guarantee that the choices of Player I
are `rectangular' in the sense that once an open set 
$U_i\subset X$ is selected, all choices $V_i\subset Y$
would then have to be combined into products $U_i\times V_i$,
and similarly for the other factor. In Plewe's game 
(\referto[Plewea], p.\ 645) (we switch the players 
to follow our original notation) this is obtained by 
letting the other player choose points $x_i\in X$,
$y_i\in Y$ in {\it alternative} steps. In our situation, 
Player I would choose, in alternative steps, 
open sets $U_i\subset X$, $V_i\subset Y$ with $x_i\in U_i$,
$y_i\in V_i$. Consider a set of choices $x_i$ by Player II 
large enough so that the corresponding sets $U_{i,x_i}$, selected 
by using a winning strategy, 
form a cover. 
Then for each such $U_{i,x_i}$, consider 
a similarly formed cover by sets of the form $V_{i,x_i,y_i}$.
The cover consisting of all rectangles of the form
$U_{i,x_i}\times V_{i,x_i,y_i}$ is in the first derivative 
of the product uniformity $\c{F}(X)\times \c{F}(Y)$. Thus,
the corresponding game $G'(X,Y,\c{V})$ is related to 
$G(X,Y,\c{V})$ in the sense that while Player I chooses 
rectangular uniform covers in the latter, the covers 
chosen in the former are uniformly locally uniform. 

In Plewe's game, the players start from an 
open cover of an open {\it rectangle}
of the product. However, as noted in his article
(\referto[Plewea], p.\ 646), for regular spaces 
this is tantamount to taking the entire product 
as the initial rectangle. Therefore, it
is now easy to see that for uniform spaces,
his game is {\it equivalent to ours with respect
to supercompleteness}. Thus, for paracompact 
factors we may state the following characterization of 
the spatiality of the localic product:

\th Theorem 6': Let $X$, $Y$ be paracompact spaces. Then 
the localic product $T(X)\times_{\loc} T(Y)$ is spatial if,
and only if, Player I has a winning strategy in the
game $G(X,Y,\c{V})$ for each open cover $\c{V}$ of 
$X\times Y$.

However, it is to be noted that the game in \referto[Plewea]
is more general than the ones described above, because 
they are not restricted to uniform spaces or paracompact
products, which always are completely regular. Nevertheless,
our characterization can be extended to 
products of general regular spaces by using the 
main result of this paper to be given in the last section. We 
obtain a deeper connection between spatiality and the locally
fine operation by moving to `covering monoids' of spaces.

\section{Formal spaces.} Motivated by locales,
Fourman and Grayson \referto[FourmanGrayson] introduced in 1981
a `formal space' of a theory, based on four conditions 
of an {\it entailment 
relation} in a propositional language, 
a pre-ordered set.
These conditions 
were taken up by Sambin in 1987 (cf.\ \referto[Sambina]) who developed 
a theory of formal spaces from a `pure' standpoint
in the spirit of the intuitionistic type theory of 
Martin-L\"of. Accordingly, intuitionistic versions 
of classical theorems for topological spaces
have been proved by several authors (see, e.\ g., Tychonoff's
Theorem in 
Negri and Valentini \referto[NegriValentini] and
Coquand's version of van der Waerden's theorem
on arithmetic progressions \referto[Coquanda].)
Formal spaces were used in 1990 by Sigstam to give 
an effective theory of spaces in her thesis \referto[Sigstama].
The approach is opposite (`top-down') to the traditional 
(`bottom-up') constructive 
approaches to say, real numbers: While the same recursive 
constructions are used, one applies them to given parts of 
a space, rather than to an assumed collection of (computable)
points.

\th Definition \ThisSection.1:
Given a pre-ordered set $(P,\leq)$, a {\it covering relation}
is a subset $\Cov \subseteq P\times 2^P$ satisfying 
the following axioms:

\vbox{\advance\leftskip by .5truecm
\item{\bf C1)} if $a\in U$, then $\Cov(a,U)$.
\item{\bf C2)} if $a\leq b$, then $\Cov(a, \{b\})$.
\item{\bf C3)} if $\Cov(a,U)$ and $\Cov(a,V)$, then 
          $\Cov(a, U\wedge V)$. Here $U\wedge V$ denotes 
          the set of elements bounded by both $U$ and $V$.          
\item{\bf C4)} if $\Cov(a,U)$ and $\Cov(u,V)$ for all $u\in U$,
          then $\Cov(a, V)$. 
}

\noindent
It is the last axiom\makeftn[[In addition to the circle of 
notions represented by 1) the locally fine operation, 2) locales
and 3) formal spaces, we may add 4) Grothendieck topologies, 
because the covering relation gives the conditions
for a Grothendieck topology on a pre-ordered set. This may be followed by
5) modal logics (see \referto[Goldblatta]), closing the 
circle with the equivalence between the modal
system S4 and the closure operation in topology, well known
since the 1930's (see \referto[McKinseya]).]]
which is directly connected with our discussion. It corresponds
to the Heyting axiom of right distributivity (characterizing 
locales) and also to the locally fine condition. Indeed, 
for a pre-uniformity $\mu$ given as a filter of coverings 
of a set $X$, we define a relation 
$R\subseteq P(X)\times P(P(X))$ by setting 
$(A,\c{U})\in R$ if there is a cover $\c{V}\in \mu$ such that 
the restriction of $\c{V}$ to $A$ refines $\c{U}$. Then 
$R$ satisfies the above conditions C1) - C3). Indeed, to see this, C1) 
is obvious because if $A$ is a member of $\c{U}$, then 
we may take the `trivial' cover $\{X\}\in \mu$ as $\c{U}$.
Condition C2) is similar, and C3) follows from the requirement 
that $\mu$ be closed under finite meets. 

On the other hand,
the {\it transitivity} 
condition C4) (called that of
{\it composition} in \referto[FourmanGrayson])
is satisfied if, and only if, the pre-uniformity
$\mu$ is locally fine, i.\ e., $\lambda\mu = \mu$.
To see this, suppose that $R$ satisfies Condition C4).
Let $\{U_i\}\in \mu$, and for each $i$, 
let $\{V^i_j\}\in \mu$. Thus, $(X,\{U_i\})\in R$, 
and for each $i$, we have $(U_i, \{U_i\cap V^i_j\})\in R$,
by the definition of $R$. By the condition under 
consideration, we obtain that $(X,\{U_i\cap V^i_j\})\in R$.
Thus, there is a member $\c{V}\in \mu$ which refines 
$\{U_i\cap V^i_j\}$, whence the latter is a member 
of $\mu$ as well. Conversely, assume that 
$\mu$ is locally fine, and suppose $(A, \c{U})\in R$, 
and for all $U\in \c{U}$, let $(U, \c{V})\in R$. 
There is $\c{U}'\in \mu$ such that 
$\c{U}'\upharpoonright A\prec \c{U}$ and for each 
$U\in \c{U}'$, there is $\c{V}_U\in \mu$ such that
$\c{V}_U\upharpoonright (U\cap A) \prec \c{V}$. The cover
$\c{W} = \cup\{\c{V}_U\upharpoonright U: U\in \c{U}'\}$ is in 
$\mu^{(1)} = \mu$, and it is easily seen that
$\c{W}\upharpoonright A \prec \c{V}$. Therefore,
$(A, \c{V})\in R$, as desired. 

The reader should note that Condition 4) above (transitivity)
is the characteristic `topological condition', expressed
in locales by the Heyting axiom and in classical topology
by the {\it idempotency} of the Kuratowski closure operator
(or by the {\it transitivity} of the corresponding relation 
between sets). In this sense, $\lambda$ corresponds to 
topology.

Given only a {\it set of generators} $G\subset P\times 2^P$,
the associated covering relation $\Cov_G$
is obtained by closing $G$ under the conditions C1) -- C4).
This means forming all Noetherian trees $T$ such that 
for each element $x$ of $T$, the immediate 
successors are derived by using one of the four 
conditions. This corresponds to the idea of using
Noetherian trees to construct `recursively defined' 
refinements of open covers of uniform spaces, in 
particular in the products of paracompact spaces.
Such constructions start from the basis of uniform 
covers, which is a commutative monoid under the operation of 
meet, and closes the collection under the 
condition of transitivity, which we have seen to be equivalent 
to the locally fine condition. By the same token, formal 
spaces are often described by giving a `formal base', 
a commutative monoid $(S, \cdot, 1)$ with unit, and 
the corresponding {\it rules of inference} equivalent 
to the above conditions C1) -- C4). For example, they could
could be given as the rules
$$1)\; {a\in U\over a\models U} \qquad
  2)\; {a\cdot b\models a} \qquad
  3)\; {a\models U\; a\models V\over a\models U\cdot V}\qquad
  4)\; {a \models U \; U\models V\over
    a\models V}. $$

We will call a pair $(P, \Cov)$ a {\it covering monoid},
if $P$ is a pre-ordered set with a unique 
maximal element $1$ and $\Cov\subset P\times 2^P$ 
is a relation closed under the conditions C1--C3.
A homomorphism between covering monoids 
$(P, \Cov)$, $(Q, \Cov')$ is a map
$f: P \to Q$ such that $(a, U)\in \Cov$ implies
$(f(a), f(U))\in \Cov'$.
With a covering monoid $(P, \Cov)$, we may associate a monoid 
$(P, \mu_{\Cov})$ of {\it covers of P} under $\Cov$, 
i.\ e., $\mu_{\Cov}$ consists of all $U\subset P$ such that
$(1,U)\in \Cov$. If $f$ -- as given above --
preserves the maximal element, 
then $f$ `restricts' to a homomorphism
$(P, \mu_{\Cov}) \to (Q, \mu_{\Cov'})$.
We denote the closure of a relation $R\subset P\times 2^P$
under C4 by $\lambda R$. This closure can be obtained 
by applying the following version of Ginsburg-Isbell 
derivation on $\Cov$: Let $\Cov^{(0)} = \Cov$, 
and given $\Cov^{(\alpha)}$, let 
$\Cov^{(\alpha+1)}$ be the collection of 
all pairs $(a, V)$ for which there is 
$(a, U)\in \Cov$ such that for all $u\in U$, 
$(u, V)\in \Cov^{(\alpha)}$. For limit ordinals
$\beta$, define $\Cov^{(\beta)} = \cup\{\Cov^{(\alpha)}:
\alpha < \beta\}$. The first stable derivative is then 
$\lambda\Cov$. This closure may also be described 
in terms of Noetherian trees: $(a,V)\in \lambda\Cov$
if, and only if, there is a Noetherian tree $T$ 
such that 1) the root of $T$ is $a$; 2) for each 
element $p$ of $T$, the immediate successors 
of $p$ form a set $U\subset P$ such that $(p,U)\in \Cov$
and 3) $V = \End(T)$.




As seen above, any pre-uniformity
$\mu$ on a set $X$ is associated with a 
covering monoid $(P(X), \Cov_\mu)$ in a 
natural way. 
Motivated by this relation, we will call pre-uniformities
{\it monoids of covers} to emphasize their formal independence
of actual pre-uniform spaces. Uniform spaces will correspond 
to {\it normal monoids} of covers $\mu$, i.\ e., in which 
for each $u\in \mu$ there is $v\in \mu$ with 
$v^2 \leq u$. Corresponding to the fine uniformity (the 
filter of all normal covers of a Tychonoff space), we 
have the {\it fine monoid of covers} of a space $X$, 
written $\c{O}(X)^*$,
consisting of all covers of $X$ with an open refinement.
This should be contrasted with the {\it fine covering 
monoid} $\c{O}(X)$ of $X$ consisting of all 
{\it pairs} $(U, \c{G})$ where $U$ is an open 
subspace of $X$ and $\c{G}$ is a cover of 
$U$ with an open refinement.
We will call a monoid of covers on a space $X$, 
written $(X, \mu)$, (super)complete if
$\lambda\mu$ is fine. In the next section, we will obtain
a product theorem which implies a far-reaching 
equivalence of locales, formal spaces and covering 
monoids (and extends our previous results on supercompleteness
to non-paracompact spaces). Let us first give two essential
lemmas on products of covering monoids. 

We note that
the product of a family $(P_i, \Cov_i)$ of covering 
monoids is a pair $(P, \Cov)$, where $P$ is the 
weak direct product of the $P_i$ consisting of 
all elements $a$ of $\Pi P_i$ with $a_i = 1_i$ for almost 
all $i$, and where $(a, U)\in \Cov$ if, and only if, 
there is for each $a_i\neq 1_i$ a pair 
$(a_i, U_i)\in \Cov_i$ such that 
$\bigwedge\{ (U_i)_i: a_i\neq 1_i\}$ refines $U$, 
where $(U_i)_i$ denotes the set of all $u\in \Pi P_j$
such that $u_i \in U_i$ and $u_j = 1_j$ for $j\neq i$.
By considering only pairs of the form 
$(1, U)$, this restricts to the usual product of 
pre-uniform spaces. Indeed, in the special 
situation in which the elements of $P_i$ are
subsets of a set $X_i$, we take the subbasis of 
$\Pi\Cov_i$ to consist of pullbacks 
$\pi_i^{-1}(a,U) = (\pi_i^{-1}[a], \pi_i^{-1}[U])$, 
where $\pi_i: \Pi X_j \to X_i$ is a projection.
In the 
general situation, we consider instead `insertions'
$q_i: P_i \to \Pi P_j$ given by $q_i(a) = (x_j)$, 
where $x_i = a$ and $x_j = 1_j$ for $j\neq i$. However,
in the following three lemmas we consider the 
(set-theoretical) situation of topological spaces.

The following Observation is obvious.

\setbox3=\hbox{\ThisSection.2\ }

\th Observation \ThisSection.2: Let $X$ be a topological space.
Then $\c{O}(F) = \c{O}(X)\upharpoonright F$ for each 
closed subpace $F\subset X$.

\th Lemma \ThisSection.3: Let $(X_i)$ be a family of topological
spaces, and let $(T(X_i), \Cov_i)$ be a corresponding family 
of covering monoids. Then $\lambda\Pi \Cov_i$ has a 
basis consisting of pairs $(a, U)$, where $U$ is a collection 
of basic open rectangles. 

\proof An inductive proof can be obtained by using the 
consecutive derivatives $\Cov^{(\alpha)}$, where 
$\Cov = \Pi \Cov_i$. The claim is clearly valid 
for $\alpha = 0$. Thus, suppose it is valid for 
$\alpha$ and let $(a, U)\in \Cov^{(\alpha + 1)}$.
Then there is a cover $V$ of $a$ such that 
$(a, V)\in \Cov$, and for each $v\in V$
a cover $W_v$ such that $(v, W_v)\in \Cov^{(\alpha)}$.
But $V$ is refined by a cover $V'$ consisting of 
open basic rectangles, and for each $v\in V$, there 
is such a refinement $W_v'$ of $W_v$. It is clear that 
the elements $v'\wedge w'$, $v'\in V'$, $w'\in W'_v$ 
form a refinement $U'$ of $U$, the elements of 
which are open basic rectangles, and $(a,U')\in \lambda\Cov$.
The case of limit ordinals is obvious.

\setbox4=\hbox{\ThisSection.4\ }
\def\LastRegularLemma{\copy4}
\th Theorem \ThisSection.4: Let $(X_i)$ be a family of regular
topological spaces. Then $(1, U)\in \lambda\Pi \c{O}(X_i)$ if,
and only if, $U\in \lambda \Pi \c{O}(X_i)^*$.

\proof We will again proceed by induction. By the definition of 
the direct product of covering monoids, $(1,U)\in \mu = \Pi \c{O}(X_i)$
if and only if $U\in \nu = \Pi \c{O}(X_i)^*$. So, suppose 
$(1,U)\in \mu^{(\alpha)}$ iff $U\in \nu^{(\alpha)}$ 
(taking $\mu, \nu$ with respect to arbitrary regular spaces). 
To show that 
this is valid for $\alpha$ replaced with $\alpha + 1$, 
it is sufficient to consider the right implication. Thus, 
let $(1,U)\in \mu^{(\alpha + 1)}$. Thus, there is $(1,V)\in \mu$
such that for each $v\in V$, we have $(v,U)\in \mu^{(\alpha)}$.
By the assumption of regularity, there is a cover
$W$ of $\Pi X_i$ by closures of basic open
rectangles in $\mu$ and hence in $\nu$ which refines $V$.
For each $w\in W$, there is an extension of $U$ to a 
cover $U_w$ of $\Pi X_i$ the restriction of which to $w$
refines $U$.
 
We may assume that $U_w\in \nu^{(\alpha)}$.
Indeed, $w$ is a (topological) product of 
regular spaces, and we may use the inductive 
hypothesis. Write $w = \bigwedge\{ \pi_i^{-1}[w_i]: i\in E\}$,
where $E$ is finite and $w_i$ is the closure of an 
open subset of $X_i$. Set $X_i' = w_i$ for $i\in E$ 
and let $X_i' = X_i$ otherwise. Then 
consider the products $\mu' = \Pi \c{O}(X_i')$,
$\nu' = \Pi \c{O}(X_i')^*$. The restriction $U'$
of $U$ to $w$ satisfies $(w, U')\in (\mu')^{(\alpha)}$,
and hence by the inductive hypothesis
$U'\in (\nu')^{(\alpha)}$. By using the 
equations $\c{O}(X_i')^* = \c{O}(X_i)^*\upharpoonright w_i$
for $i\in E$, and recalling that the 
Ginsburg-Isbell derivatives preserve substructures
(i.\ e., $(\xi\upharpoonright A)^{(\alpha)} = 
  (\xi)^{(\alpha)}\upharpoonright A$)
it easily follows there is 
cover $U_w\in \nu^{(\alpha)}$ the restriction of 
which to $w$ refines $U'$, as desired.

But then the elements $w\wedge u_w$, $u_w\in U_w$, form a
cover $U''$ such that $U''\in \nu^{(\alpha + 1)}$ and 
$U''$ refines $U$, implying $U\in \nu^{(\alpha + 1)}$.
As above, the limit ordinal case is obvious.

\section{A general product theorem.}
It can be seen from the previous section 
that the theory of formal spaces corresponds to 
that of locally fine covering monoids. In this section, 
we will use notions and lemmas developed above
to link supercompleness in products
to spatiality in a general fashion. We extend the 
characterization of supercompleteness in
a paracompact product $\c{F}X\times \c{F}Y$ by the equation
$$ \lambda(\c{F}X\times \c{F}Y) = \c{F}(X\times Y)$$
to a similar one (\ThisSection.6) characterizing spatiality, even without 
paracompactness. 

We will first describe the locale 
product simply as the locally fine (or $\lambda$-)
product. The product theorem given in this section grew 
out of the author's attempt to understand the proof
given by Dowker and Strauss (\referto[DowkerStrauss])
for their product theorem. The following  definitions 
are well-known, see, e.\ g., \referto[Sigstama].

Let $\Cov$ be a covering relation on a pre-ordered set $P$.
For subsets $U,V \subset P$, define $U \leq V$ if
for all $u\in U$, we have $\Cov(u, V)$. Then define 
an equivalence relation $\sim$ by setting $U\sim V$
if $U\leq V$ and $V\leq U$. Denote the 
corresponding equivalence classes by $[U]$. The locale associated 
with the covering relation $\Cov$ is the set 
$$ \c{L}_{\Cov} \; = \; \{[U]: U\subseteq P\}$$
equipped with the lattice operations (recall the definition
of $U\wedge V$)
$$ [U] \wedge [V] \; = \; [U\wedge V],\;
\bigvee_{i\in I}[U_i] \; = \; [\cup\{U_i: i\in I\}].$$
We say that the covering relation $\Cov$ 
(or more exactly the pair $(P, \Cov)$) {\it generates $L$}.
We extend this definition to covering monoids
by stipulating that the locale associated with a covering 
monoid $(P, \mu)$ is the one generated by $\lambda\mu$.

On the other hand, given a locale $L$, define a canonical 
covering relation $\Cov_L$ by setting 
$\Cov_L(a,U) \Leftrightarrow a\leq \bigvee U$. Then 
$\c{L}_{Cov_L} \cong L$. Thus, every locale has 
a canonical generating covering relation, and it follows 
from the right distributivity of the locale that 
this relation is locally fine, i.\ e., defines a 
formal space. If $(P, \Cov)$ generates $L$, then 
there is a canonical embedding (of covering monoids)
$(P,\Cov) \hookrightarrow (2^P, \Cov_L)$ given 
by $a\mapsto [a]$; we will consider the generating 
monoid a submonoid of $(2^P, \Cov_L)$.

One says a subset $U\subseteq L$ is a {\it cover} of a 
locale $L$ if $\bigvee U = 1$. A subset $V$ is a 
{\it refinement} of $U$ if for each $v\in V$ there is 
$u\in U$ such that $u \leq v$. We denote the monoid 
of all covers of $L$ by $\Cov(L)$. Thus, $\Cov(L)$ 
is the collection of all $U\subseteq L$ such that 
$\Cov_L(L,U)$, and by transitivity, $\Cov(L)$ is 
locally fine. (Note that since $1\in L$, 
$\Cov_L(L,U)$ implies $\bigvee U = 1$.) 

We will construct
the co-product $\amalg L_i$ of 
given frames $L_i$. (We remind the reader that 
the difference with locales is that morphisms 
go in the opposite direction. With the 
product of locales, we have projections 
$\pi_j: \Pi L_i\to L_j$, whereas with the co-product
of frames, we have `insertions' $q_j: L_j\to \amalg L_i$.)

Let $(L_i)$ be a family of frames, and consider the 
Cartesian product frame $\Pi L_i$. (This is a frame, 
but not the {\it product} of the $L_i$ in the 
category of locales!) 
Take a subframe $B\subset \Pi L_i$
which consists of all $b = (b_i)$ such that 
$b_i = 1_i$ except for finitely many $i$ (the direct product 
of monoids). 
We define a covering relation by first choosing a 
set $G\subset (B, 2^B)$ of generators to 
consist 
of all $(a,U)$, where for some $i$, $U$ has the form:
For $j\neq i$, $\pi_j[U] = \{a_j\}$, and 
$(a_i, \pi_i[U])\in \Cov_{L_i}$. Thus, $U$ has been 
obtained from $a$ by `splitting' it along
exactly one coordinate direction with respect to 
the corresponding covering relation. (Notice that this 
condition corresponds to the 
`coordinatewise derivation condition' from Section 2.)
The associated 
covering relation $\Cov_G$ is obtained by closing 
$G$ under the conditions C1)--C4) of covering relations.
The frame $L = \c{L}_{\Cov_G}$ will be our co-product.
Recall that the elements of $L$ are equivalence 
classes $[U]$ of subsets $U\subset B$ under the 
equivalence relation: $U\sim V$ iff $U\leq V$ and
$V\leq U$. It follows that $U$ is a cover of 
$L$, i.\ e., $\bigvee U = 1_L$, if 
$(1,U)\in \Cov_G$, where 1 denotes the maximal 
element of $B$. We will show that for each such 
$U$ there is $\c{U}\in \lambda\Pi\mu_i$ with 
$\phi(\c{U}) \prec U$, where $\mu_i = \Cov(L_i)$ and 
$\phi: Cov(\lambda\Pi\mu_i)\to Cov(L)$ is an
embedding of covering monoids. 

The product of the $\mu_i$ has a subbasis defined 
by the insertions $q_i: L_i \to B$ by setting
$ q_i(x) = (a_j)$, where $a_i = x$ and $a_j = 1_j$ for 
$j\neq i$. It is obvious that $(a, U)\in Cov_{L_i}$
implies $(q_i(a),q_i(U))\in G$, by the definition of $G$.
By taking finite meets, it turns out that 
$\Pi \mu_i$ has a basis $\c{B}$ contained in $\Cov_G$. 
We note that $\Pi\mu_i$ is obtained from $\c{B}$ by 
applying the rules C1--C3.

How are the elements $(a,U)\in \Cov_G$ obtained?
By the definition, $(a,U)$ belongs to $\Cov_G$ if, 
and only if, there is Noetherian tree $T$ such that 
1) the root of $T$ is a; 2) the immediate 
successors of an element $p\in T$ are obtained from 
$p$ by applying $G$ or one of the conditions C1--C4, 
and 3) $U =\End(T)$. It is clear that 
$\lambda\Pi \mu_i$ is closed under these conditions.
Thus, it is contained in $\Cov_G$. The opposite inclusion 
is clear, too, and hence $\Cov_G$ and $\Pi\mu_i$ are the 
same covering relation.
  
Let $\c{E}$ be a cover of $L$. Thus, $\bigvee \c{E} = 1_L$.
Hence, there are sets $E_i \subset B$ such that 
$\c{E} = \{[E_i]: i\in I\}$, and therefore 
$$ 1_L \;=\; \bigvee_{i\in I}[E_i] \; = \; 
                        \left[\cup_{i\in I} E_i\right].$$
It follows that 
$ (1, \cup_{i\in I}E_i)\in \Cov_G,$
where 1 denotes the maximal element of $B$. Thus, 
$ \c{U} = \cup\{E_i: i\in I\}$ is an element of
$\lambda\Pi R_i$ such that $[\c{U}] = \{[u]: u\in \c{U}\}$
refines $\c{E}$. Denote the covering monoid associated 
with $\Cov_G$ by $\Cov(G)$. The mapping 
$u \mapsto [u]$ defines a natural homomorphism 
$\phi: \Cov(\Pi\mu_i) \to \Cov(\amalg L_i)$ 
of covering monoids. However, the factors $L_i$
are frames and hence partially ordered and so 
is the weak direct product $B$. It follows that 
$u\mapsto [u]$ yields an embedding $B \to L$, 
which extends to covers. Thus, we have proved the 
following result:

\th Theorem \ThisSection.1: Let $(L_i)$ be a family 
of frames. Then there is an embedding $\phi$ of 
covering monoids
$$ \lambda(\Pi \Cov_{L_i}) \to^{\phi} \Cov_{\amalg L_i}\leqno(*) $$
where $\Cov_{L_i}$ is the canonical covering relation 
on $L_i$ and the left-hand side of $(*)$ is a locally fine 
generating covering monoid for $\amalg L_i$.
Moreover, the mapping $\phi$ is induced by the 
embedding $u\mapsto [u]$, and 
for each cover $\c{V}$ of 
$\amalg L_i$ there is $\c{U}\in \lambda\Pi \Cov_{L_i}$
such that $\phi(\c{U}) \prec \c{V}$. 

For pre-uniform spaces $(X_i, \mu_i)$, 
the direct product is a pair $(\Pi X_i, \Pi \mu_i)$, 
where $\Pi \mu_i$ is generated by the basis of all 
finite meets of single pullbacks $\pi^{-1}_i[\c{U}]$,
$\c{U}\in \mu_i$ (`basic rectangular covers').
Moreover, $\lambda\Pi \mu_i$ is generated by 
covers consisting of basic rectangles, which
form a monoid. Lacking better notation, we denote
this monoid of rectangles by $[\lambda\Pi \mu_i]_{\c{R}}$.
Its covering relation is induced by the 
pre-uniform structure of the product: A collection of 
rectangles cover a rectangle if, and only if, they cover
the latter as a (pre-)uniform cover.

This is 
special case of the product of covering monoids 
$(P_i, \mu_i)$, in which the basic rectangular 
covers are finite meets of pullbacks of the form
$\pi_i^{-1}(a,U) = (\pi_i^{-1}[a], \pi_i^{-1}[U])$, where 
$\pi_i^{-1}[a]$ is a basic open rectangle covered 
by the cover $\pi_i^{-1}[U]$ consisting of 
basic rectangles. 
\ThisSection.1 
implies the following result.

\th Corollary \ThisSection.2: Let $(X_i)$ be a 
family of topological spaces. For each $X_i$,
let $\c{O}_i(X_i)$ denote the fine covering monoid 
induced by open-refinable covers.
Then there is an embedding $\phi$
$$ [\lambda\Pi \c{O}(X_i)]_{\c{R}} \; \cong \; 
      \lambda\Pi \Cov_{T(X_i)}\;\hookrightarrow^{\phi}\;
           \Cov_{\amalg T(X_i)},$$
where $T(X_i)$ is the topology of $X_i$. Moreover, 
for any cover $\c{V}$ of $\amalg T(X_i)$, there is 
a rectangular cover $\c{U}$ in $\lambda\Pi \c{O}_i(X_i)$
such that $\phi(\c{U}) \prec \c{V}$.

Since $\lambda\Pi \Cov_{T(X_i)}$ generates the localic 
product of the $T(X_i)$, we (ab)use the above corollary 
to say that $\lambda\Pi \c{O}(X_i)$ generates
it, too.

In \referto[Isbellc], Isbell showed that the product 
of paracompact locales is paracompact. 
Dowker and Strauss \referto[DowkerStrauss]
extended this result to include the cases of metacompact 
and Lindel\"of (regular) locales. These results (and 
an unlimited number of others) follow from Theorem
\ThisSection.2. 

Indeed, for a topological space $X$, 
let $\c{LF}(X)$ be the monoid of all covers which have an 
open, locally finite refinement. Let us call $\c{LF}(X)$
the locally finite monoid of covers on $X$. Then 
$\c{LF}(X)$ is locally fine, and $X$ is paracompact if
$\c{LF}(X)$ is fine, i.\ e., contains (and thus equals)
$\c{O}(X)$. It is easy to see that arbitrary direct products 
of locally finite monoids of covers is again locally finite.
(We call $\mu$ on a space $X$ locally finite if it contains
$\c{LF}(X)$.) This follows from the easy observation 
that any binary, and more generally finite, product of 
locally finite covers is again locally finite. Finally, 
$\lambda$ preserves local finiteness, so that 
$\lambda\Pi \c{LF}(X_i)$ is locally finite. The same is 
true of point-finite, locally countable, point-countable,
Lindel\"of, and compact monoids of covers (call $\mu$
compact if every cover has a finite open refinement in 
$\mu$). These considerations are valid 
for general covering monoids. 
Therefore, we obtain the following corollary:

\th Corollary \ThisSection.3: If the members of a family 
$(L_i)$ of locales are
compact (resp.\ paracompact, metacompact,
Lindel\"of, para-Lindel\"of, meta-Lindel\"of), then so 
is $\Pi_{\loc} L_i$.

In fact, we may use \ThisSection.2 to establish a 
relation between 
the spatiality of localic products and the locally fine
condition. To this end, we might first give a game-theoretical
characterization for the $\lambda$ of the product of fine 
monoids to be fine, and show its equivalence with Plewe's 
game-theoretical characterization of spatiality in products.
However, we will proceed directly. Let $(X_i)$ be a family 
of sober spaces, i.\ e., $\Pt(T(X_i) \cong X_i$. We will first
show that the localic product $\Pi_{\loc} T(X_i)$ is spatial,
$\Pi_{\loc} T(X_i) = T(\Pi X_i)$ if, and only if, 
$\lambda\Pi\c{O}(X_i)$ is the fine monoid $\c{O}(\Pi X_i)$.

\th Theorem \ThisSection.4: The localic product of 
a family $(X_i)$ of sober topological spaces is spatial
if, and only if, $\lambda\Pi \c{O}(X_i) = \c{O}(\Pi X_i)$.

\proof Suppose that 
$\lambda\Pi\c{O}(X_i) = \c{O}(\Pi X_i)$.
The locales $T(X_i)$ are generated by the fine 
covering monoids $\c{O}(X_i)$. Hence, their 
localic product $\Pi_{\loc} T(X_i)$ is generated by 
$\lambda\Pi \c{O}(X_i)$, which is by assumption the 
fine monoid of the topological product, and hence 
generates $T(\Pi X_i)$, as desired. On the other
hand, suppose that $\Pi_{\loc} T(X_i)$ is spatial.

Given an open cover $U$ of a basic open rectangle $a$ in $\Pi X_i$, 
we may consider $U$ a cover of $a$ in 
$L = \Pi_{\loc} T(X_i)$. 
But $L$ is generated by $\lambda\Pi \c{O}(X_i)$, 
and hence there is, for each $u\in U$, collection $V_u$ of open sets
(basic rectangles) such that $u = [V_u]$ and 
hence 
$(a, V)\in \lambda\Pi \c{O}(X_i)$, where 
$V = \cup\{V_u: u\in U\}$. Therefore, $V$ is a 
refinement of $U$ in the locally fine closure 
of the product of the $\c{O}(X_i)$, which 
consequently  refines the fine monoid 
of the topological product, i.\ e., it is itself fine.

Notice in particular that we have not assumed the factors
to be regular. However, this result cannot be directly 
applied to spaces (via spatiality) along the lines 
of \ThisSection.3, because the fine covering monoids 
$\c{O}(X)$ carry -- within their structure -- all
the open subspaces. As a consequence, after taking 
the locally fine coreflection, the corresponding 
products $\Pi\c{O}(X_i)$ produce in general monoids of covers 
finer than the ones obtained from products of 
monoids $\c{O}(X_i)^*$ of covers (as generalized 
pre-uniform spaces). In order to bridge the gap, 
we need to assume regularity. The following 
lemma provides a link between spatiality and 
$\lambda$-covers.

\th Lemma \ThisSection.5: Let $(X_i)$ be a family 
of topological spaces, and let $\c{U}$ be a 
collection of basic open rectangles in $\Pi X_i$
such that $\c{U}' = \{[u]: u\in \c{U}\}$ covers
the points of $\Pi_{\loc} T(X_i)$. If $\c{U}$ belongs to 
$\lambda\Pi\c{O}(X_i)^*$, then $\c{U}'$ covers $\Pi_{\loc} T(X_i)$.

\proof This follows immediately from the result that 
$\lambda\Pi\c{O}(X_i)$ generates $\Pi_{\loc} T(X_i)$.

The condition that 
$\lambda\Pi \c{O}(X_i)^* = \c{O}(\Pi X_i)^*$ is 
analogous to the condition -- studied by the author -- that 
$\lambda(\Pi\c{F}(X_i))$ contain all normal
covers of $\Pi X_i$. (In \referto[Pelanta]
this was shown to be true whenever the $X_i$
are completely metrizable spaces; in \referto[HohtiHusekPelant],
this result has been extended to paracompact spaces which are 
countable unions of closed, partition-complete
subspaces.)

\th Theorem \ThisSection.6: The localic product of 
a family $(X_i)$ of regular topological spaces is spatial
if, and only if, $\lambda(\Pi \c{O}(X_i)^*) = \c{O}(\Pi X_i)^*$.

\proof Suppose that $\Pi_{\loc} T(X_i)$ is spatial. Then 
by \ThisSection.4, 
$\lambda\Pi\c{O}(X_i) = \c{O}(\Pi X_i)$, and hence
by \LastRegularLemma we have 
$\lambda(\Pi \c{O}(X_i)^*) = \c{O}(\Pi X_i)^*$, as required.

On the other hand, suppose that this condition holds.
We recall that a regular locale is spatial iff it 
does not contain a non-empty, closed pointless
sublocale. If $\Pi_{\loc} T(X_i)$ were not spatial,
then there would exist such a sublocale $F$, and 
the collection $\c{U}$ of all basic open rectangles $u$
of $\Pi X_i$ such that $[u]\wedge F = 0$ would form 
an open cover $\Pi X_i$ for which $\c{U}' = 
\{[u]: u\in \c{U}\}$ covers the points of the 
localic product. But by the assumption 
$\c{U}\in \lambda\Pi\c{O}(X_i)^*$, and hence (by the preceding lemma)
$\c{U}'$ would cover the product locale, which
is impossible. Hence, the product in question is 
spatial.

\bigskip
\SectionBreak
\centerline{\smallcaps References}
\bigskip
{
\beginref

\ref Alster, K: A class of spaces whose Cartesian product with
     every hereditarily Lindel\"of space is Lindel\"of.-
     Fund.\ Math.\ 114:3, 1981, pp.\ 173--181 [Alstera]

\ref Arhangel'skii, A: On topological spaces complete in the sense
     of \v Cech.- Vestnik Moskov.\ Univ.\ Ser.\ I.\ Mat.\ Mekh.\
     2, 1961, pp.\ 37--40 (Russian) [Arhangelskiia]

\ref Benabou, J: Treillis locaux et paratopologies.- S\'eminaire
     C.\ Ehresmann de Topologie et de G\'eometrie
     Diff\'e\-rentiel\-le, 1957/58, Fac.\ des Sciences de Paris, 1959
     [Benaboua]

\ref Coquand, T: Minimal invariant spaces in formal topology.-
     J.\ Symb.\ Logic 62:3, 1997, pp.\ 689--698
     [Coquanda]

\ref Corson, H: Determination of paracompactness by uniformities.-
     Amer.\ J.\ Math. 80, 1958, pp.\ 185--190 [Corsona]

\ref Dowker, C.\ H., and D.\ Strauss: Sums in the category
     of frames.- Houston J.\ Math. 3, 1977, pp.\ 17--32 [DowkerStrauss]

\ref Ehresmann, C: Gattungen von lokalen strukturen.- Jahresbericht
     Deutsch.\ Math.\ Verein.\ 60, 1957, pp.\ 59--77 
     [Ehresmanna]

\ref Fletcher, P., and W.\ F.\ Lindgren: $C$-complete quasi-uniform
     spaces.- Arch.\ Math.\ (Basel) 30:2, 1978, pp.\ 175--180
     [FletcherLindgren]

\ref Fourman, M., and R.\ Grayson: Formal spaces.- The L.\ E.\ J.\
     Brouwer Centenary Symposium, A.\ S.\ Troelstra and D.\ van Dalen
     (eds.), North-Holland, 1982, pp.\ 107--122 [FourmanGrayson]

\ref Friedler, L.\ M., H.\ W.\ Martin and S.\ W.\ Williams:
     Paracompact C-scattered spaces.- Pacific Math.\ Journal 129:2,
     1987, pp.\ 277 -- 296 [Friedlera]

\ref Frol\'\i k, Z: On the topological product of paracompact spaces.-
     Bull.\ Acad.\ Pol.\ Sci.\ Math.\ 8, 1960, pp.\ 747 -- 750
     [Frolika]

\ref Ginsburg, S.\ and J.\ R.\ Isbell: Some operators on uniform
     spaces.- Trans.\ Amer.\ Math.\ Soc.\ 93, 1959, pp.\ 145 -- 168
     [Ginsburga]

\ref Goldblatt, R: Grothendieck topology as geometric modality.- 
     Z.\ Math.\ Logik Grundlag.\ Math.\ 27:6, 1981, pp.\ 495--529
     [Goldblatta]

\ref Hager, A.\ W: Some nearly fine uniform spaces.- Proc.\ London
     Math.\ Soc.\ (3), 28, 1974, pp.\ 517--546 [Hagera]

\ref Hausdorff, F: Die Schlichten stetigen Bilden des Nullraums.-
     Fund.\ Math.\ XXIX, 1937, pp.\ 151--158 [Hausdorffa]

\ref Hohti, A: On uniform paracompactness.- Ann.\ Acad.\ Scient.\
     Fenn., Series A, I.\ Mathematica, Dissertationes 36, 1981
     [Hohtia]

\ref Hohti, A: On supercomplete uniform spaces (II).- Reports
     of the Department of Mathematics, University of Helsinki,
     1982  [HohtiReport]

\ref Hohti, A:On supercomplete uniform spaces.- Proc.\ Amer.\ Math.\ Soc.\ 87:,
     1983, pp.\ 557--560 [Hohtib]

\ref Hohti, A: Uniform hyperspaces.- Colloquia Mathematica Societatis
     J\'anos Bolyai 41, Topology and Applications (Eger, Hungary, 1983),
     pp.\ 333--344 [HohtiEger]

\ref Hohti, A: On Ginsburg-Isbell derivatives and ranks of metric spaces.-
     Pacific J.\ Math.\ 111 (1), 1984, pp.\ 39--48 [Hohtic]

\ref Hohti, A: On supercomplete uniform spaces II.- Czechosl.\ Math.\ J.\ 37,
     1987, pp.\ 376--385 [Hohtid]

\ref Hohti, A: On supercomplete uniform spaces III.- Proc.\ Amer.\ Math. Soc.\
     97:2, pp.\ 339--342 [Hohtie]

\ref Hohti, A: On Relative $\omega$--cardinality and Locally Fine Coreflections
        of Products.- Topology Proceedings 13:1, 1989, pp.\ 93--106 [Hohtif]

\ref Hohti, A: On supercomplete spaces V: Tamano's product problem.-
     Fund.\ Math.\ 136:2, 1990, pp.\ 121--125 [Hohtig]

\ref Hohti, A., and J.\ Pelant: On complexity of metric
     spaces.- Fund.\ Math.\
     CXXV, 1985, pp.\ 133--142 [HohtiPelanta]

\ref Hohti, A., and J.\ Pelant: On supercomplete uniform spaces IV:
     countable products.-  Fund.\ Math. 136:2, 1990, pp.\ 115--120
     [HohtiPelantb]

\ref Hohti, A., and Yun Z: Countable products of \v Cech-scattered supercomplete
     spaces.- Czechosl.\ Math.\ J., to appear [HohtiYun]

\ref Hohti, A., Hu\v sek, M., and J.\ Pelant: The Locally Fine
     Coreflection and Normal Covers in the Products of
     Partition-complete Spaces.- manuscript. [HohtiHusekPelant]

\ref Hu\v sek, M., and J.\ Pelant: Extensions and restrictions
     in products of metric spaces.- Topology Appl.\ 25, 1987,
     pp.\ 245 -- 252 [HusekPelant]

\ref Isbell, J: Supercomplete spaces.- Pacific J.\ Math.\ 12, 1962,
     pp.\ 287 -- 290 [Isbella]

\ref Isbell, J: Uniform spaces.- Math.\ Surveys, no.\ 12, Amer.\
     Math.\ Soc., Providence, R.\ I., 1964 [Isbellb]

\ref Isbell, J: Atomless parts of spaces.- Math.\ Scand.\ 31, 1972,
     pp.\ 5--32 [Isbellc]

\ref Isbell, J: Product spaces in locales.- Proc.\ Amer.\ Math.\ Soc.,
     81:1, 1981, pp.\ 116--118 [Isbelld]

\ref Johnstone, P.\ T: Stone spaces.- Cambridge Univ.\ Press, 
     1982 [Johnstonea]

\ref Johnstone, P.\ T: The point of pointless topology.-
     Bull.\ Amer.\ Math.\ Soc.\ 8:1, 1983, pp.\ 41--53
     [Johnstoneb]

\ref McKinsey, J.\ C.\ C: A solution of the decision problem
     for the Lewis systems S2 and S4 with an application to 
     to topology.- J.\ Symb.\ Logic 6, 1941, pp.\ 117--134
     [McKinseya]

\ref Michael, E: A note on completely metrizable spaces.-
     Proc.\ Amer.\ Math.\ Soc.\ 96:3, 1986, pp.\ 513--522 [Michaelb]

\ref Negri, S., and S.\ Valentini: Tychonoff's theorem in the
     framework of formal topologies.- J.\ Symb.\ Logic 62:4, 1997,
     pp.\ 1315 -- 1332 [NegriValentini]

\ref Papert, D., and S.\ Papert: Sur les treillis des ouverts
     et les paratopologies.- S\'eminaire
     C.\ Ehresmann de Topologie et de G\'eometrie
     Diff\'erentielle, 1957/58, Fac.\ des Sciences de Paris, 1959
     [Paperta]

\ref Pelant, J: Locally fine uniformities and normal covers.-
     Czechosl.\ Math.\ J.\ 37 (112), 1987, pp.\ 181 -- 187
     [Pelanta]

\ref Plewe, T: Localic products of spaces.- Proc.\ London Math.\
     Soc.\ 73:3, 1996, pp.\ 642--678 [Plewea]

\ref Plewe, T:Countable products of absolute $C_\delta$-spaces.-
     Topology Appl.\ 74, 1996, pp.\ 39--50 [Pleweb]

\ref Przymusinski, T: On the notion of $n$-cardinality.-
     Proc.\ Amer.\ Math.\ Soc.\ 69, 1978, pp.\ 333--338
     [Przymusinskia]

\ref Rudin, M.\ E., and S.\ Watson: Countable products of
     scattered paracompact spaces.- Proc.\ Amer.\ Math.\ Soc.\
     89:3, 1983, pp.\ 551 -- 552 [RudinWatson]

\ref Sambin, G: Intuitionistic formal spaces and their
     neighbourhood.- Logic Colloquium '88, R.\ Ferro et al
     (eds.), North-Holland, Amsterdam, 1989, pp.\ 261--285 [Sambina]

\ref Sigstam, I: Formal spaces and their effective presentations.-
     Arch.\ Math.\ Logic 34:4, 1995 pp.\ 211 -- 246 [Sigstama]

\ref Stone, A.\ H: Kernel constructions and Borel sets.- Trans.\
     Amer.\ Math.\ Soc.\ 107, 1963, pp.\ 58--70 [Stonea]

\ref Telg\'arsky, R: C-scattered and paracompact spaces.- Fund.\
     Math.\ 73, 1971, pp.\ 59 -- 74 [Telgarskya]

\ref Telg\'arsky, R., and H.\ H.\ Wicke: Complete exhaustive sieves
     and games.- Proc.\ Amer.\ Math.\ Soc.\ 102:1, 1988, pp.\
     737--744 [TelgarskyWicke]

\endref
}

\vskip2cm
\noindent
{\bf The address:}
\bigskip
{\ninepoint\obeylines
\hskip\refskip University of Helsinki
\hskip\refskip Department of Mathematics
\hskip\refskip PL 4 (Yliopistonkatu 5)
\hskip\refskip 00014 HELSINGIN YLIOPISTO
\hskip\refskip FINLAND
}

\end